\documentclass[11pt,a4paper,reqno]{amsart}
\usepackage[T1]{fontenc}
\usepackage{color,latexsym,amsfonts,amssymb,bbm,comment}
\usepackage{hyperref}
\usepackage{amsmath,cite}
\usepackage{amsthm}
\usepackage{fullpage}
\usepackage{graphicx}
\usepackage{tikz}
\usepackage{caption,subcaption}
\usepackage{makecell}
\usepackage{boldline}
\setcellgapes{3pt}
\usepackage{multirow}
\usepackage{mathtools}
\usepackage{mathrsfs}

\newtheorem {thm}{Theorem}[section]
\newtheorem {prop}[thm]{Proposition} 

\newtheorem {lem}[thm]{Lemma}

\def\N{{\Bbb N}}

\def\Z{{\Bbb Z}}
\def\R{{\Bbb R}}
\def\one{\mathbbmss{1}}

\def\P{{\Bbb P}}

\def\E{{\Bbb E}}

\def\e{{\varepsilon}}

\def\diam{{\rm diam}}

\def\eps{\varepsilon}
\def\dist{{\rm dist}}

\def\a{\alpha}

\def\d{\delta}

\def\e{\varepsilon}

\def\phi{\varphi}

\def\g{\gamma}

\def\k{\kappa}

\def\r{\rho}

\def\x{\xi}

\def\G{\Gamma}

\def\P{{\Phi}}

\def\T{\T}

\def\P{{\mathcal P}}

\def\V|{{\Vert}}

\def\laI{\la_{\rm I}}

\def\laW{\la_{\rm W}}
\def\muW{\mu_{\rm W}}
\def\muS{\mu_{\rm S}}
\def\XS{X_{\rm S}}
\def\XW{X_{\rm W}}
\def\muR{\mu_{r}}
\def\Pspa{\mathcal{P}^o}
\def\Espa{\mathcal{E}^o}
\def\Pinf{P^{\xi}}
\def\Einf{E^{\xi}}
\def\EExt{E}
\def\ELoc{L}
\def\EGlob{G}
\def\ESurv{E^{\rm{c}}}

\pgfmathsetseed{14}

\newcommand{\coxSegment}[6]{

	\pgfmathsetmacro{\Xa}{#1}
	\pgfmathsetmacro{\Xb}{#2}
	\pgfmathsetmacro{\Xc}{#3}
	\pgfmathsetmacro{\Xd}{#4}

	\draw (\Xa - \factor*\Xc + \factor*\Xa,\Xb  - \factor*\Xd + \factor*\Xb) -- (\Xc + \factor*\Xc - \factor*\Xa, \Xd + \factor*\Xd - \factor*\Xb);
	\foreach \y in {1,2,...,#5}{
		\pgfmathsetmacro{\Xf}{random()}
		\pgfmathsetmacro{\Xg}{1- \Xf)}
		\pgfmathsetmacro{\Xh}{\Xa*\Xf + \Xc*\Xg}
		\pgfmathsetmacro{\Xi}{\Xb*\Xf+ \Xd*\Xg}
		\pgfmathparse{\Xh*\Xh + \Xi*\Xi > #6 ? 1: 0}
		\ifthenelse{\pgfmathresult>0}{\fill (\Xh, \Xi) circle (2pt);}{}
	}

}

\newcommand{\coxSegmentBoxes}[3]{
	\foreach \y in {1,2,...,#1}{
		\pgfmathsetmacro{\Yf}{random()}
		\pgfmathsetmacro{\Yg}{random()}
		\pgfmathparse{\Yf*\Yf + \Yg*\Yg > #2/(#3*#3) ? 1: 0}
		\ifthenelse{\pgfmathresult>0}{\fill[gray] (#3*\Yf, #3*\Yg) circle (2pt);}{}	
	}
		\foreach \y in {1,2,...,#1}{
		\pgfmathsetmacro{\Yf}{random()}
		\pgfmathsetmacro{\Yg}{-random()}
		\pgfmathparse{\Yf*\Yf + \Yg*\Yg > #2/(#3*#3) ? 1: 0}
		\ifthenelse{\pgfmathresult>0}{\fill[gray] (#3*\Yf, #3*\Yg) circle (2pt);}{}	
	}
		\foreach \y in {1,2,...,#1}{
		\pgfmathsetmacro{\Yf}{-random()}
		\pgfmathsetmacro{\Yg}{random()}
		\pgfmathparse{\Yf*\Yf + \Yg*\Yg > #2/(#3*#3) ? 1: 0}
		\ifthenelse{\pgfmathresult>0}{\fill[gray] (#3*\Yf, #3*\Yg) circle (2pt);}{}	
	}
		\foreach \y in {1,2,...,#1}{
		\pgfmathsetmacro{\Yf}{-random()}
		\pgfmathsetmacro{\Yg}{-random()}
		\pgfmathparse{\Yf*\Yf + \Yg*\Yg > #2/(#3*#3) ? 1: 0}
		\ifthenelse{\pgfmathresult>0}{\fill[gray] (#3*\Yf, #3*\Yg) circle (2pt);}{}	
	}

}

\newcommand{\coxSegmentBoxesG}[3]{
	\foreach \y in {1,2,...,#1}{
		\pgfmathsetmacro{\Yf}{random()}
		\pgfmathsetmacro{\Yg}{random()}
		\pgfmathparse{\Yf*\Yf + \Yg*\Yg > #2/(#3*#3) ? 1: 0}
		\ifthenelse{\pgfmathresult>0}{\pgcircle{#3*\Yf}{#3*\Yg}{0.2cm};}{}	
	}
		\foreach \y in {1,2,...,#1}{
		\pgfmathsetmacro{\Yf}{random()}
		\pgfmathsetmacro{\Yg}{-random()}
		\pgfmathparse{\Yf*\Yf + \Yg*\Yg > #2/(#3*#3) ? 1: 0}
		\ifthenelse{\pgfmathresult>0}{\pgcircle{#3*\Yf}{#3*\Yg}{0.2cm};}{}	
	}
		\foreach \y in {1,2,...,#1}{
		\pgfmathsetmacro{\Yf}{-random()}
		\pgfmathsetmacro{\Yg}{random()}
		\pgfmathparse{\Yf*\Yf + \Yg*\Yg > #2/(#3*#3) ? 1: 0}
		\ifthenelse{\pgfmathresult>0}{\pgcircle{#3*\Yf}{#3*\Yg}{0.2cm};}{}	
	}
		\foreach \y in {1,2,...,#1}{
		\pgfmathsetmacro{\Yf}{-random()}
		\pgfmathsetmacro{\Yg}{-random()}
		\pgfmathparse{\Yf*\Yf + \Yg*\Yg > #2/(#3*#3) ? 1: 0}
		\ifthenelse{\pgfmathresult>0}{\pgcircle{#3*\Yf}{#3*\Yg}{0.2cm};}{}	
	}

}

\def\eq{\begin{equation}}
\def\en{\end{equation}}
    \def\e{{\varepsilon}}

    \def\eps{\varepsilon}
    \def\dist{{\rm dist}}

    \def\a{\alpha}

    \def\e{\varepsilon}
    
    \def\phi{\varphi}
    \def\g{\gamma}
    \def\la{\lambda}
    \def\k{\kappa}
    
    \def\r{\rho}

    \def\x{\xi}

    \def\G{\Gamma}
    
    \def\P{{\Phi}}
    
    \def\T{\T}

    \def\V|{{\Vert}}

    \def\d{{\rm d}}
    \def\E{\mathbb{E}}
    \def\V{\mathbb{V}}

    \def\one{\mathbbmss{1}}

    \def\N{\mathbb{N}}
    \def\one{\mathbbmss{1}}
    \def\P{\mathbb{P}}
    \def\R{\mathbb{R}}
    \def\Z{\mathbb{Z}}

\keywords{Interacting particle systems; random graphs; survival; extinction; percolation; Boolean model}
\subjclass[2010]{Primary 60J25; secondary 60K35, 60K37}

\begin{document}
\author{Alexander Hinsen}
\address[Alexander Hinsen]{Weierstrass Institute Berlin, Mohrenstr. 39, 10117 Berlin, Germany}
\email{Alexander.Hinsen@wias-berlin.de}

\author{Benedikt Jahnel}
\address[Benedikt Jahnel]{Weierstrass Institute Berlin, Mohrenstr. 39, 10117 Berlin, Germany}
\email{Benedikt.Jahnel@wias-berlin.de}

\author{Elie Cali}
\address[Elie Cali]{Orange SA, 44 Avenue de la R\'epublique, 92326 Ch\^atillon, France}
\email{Elie.Cali@orange.com} 

\author{Jean-Philippe Wary}
\address[Jean-Philippe Wary]{Orange SA, 44 Avenue de la R\'epublique, 92326 Ch\^atillon, France}
\email{Jeanphilippe.Wary@orange.com}

\title{Phase transitions for chase-escape models\\ on Gilbert graphs}

\date{\today}

\maketitle

\begin{abstract}
We present results on phase transitions of local and global survival in a two-species model on Gilbert graphs. At initial time there is an infection at the origin that propagates on the Gilbert graph according to a continuous-time nearest-neighbor interacting particle system. The Gilbert graph consists of susceptible nodes and nodes of a second type, which we call white knights. The infection can spread on susceptible nodes without restriction. If the infection reaches a white knight, this white knight starts to spread on the set of infected nodes according to the same mechanism, with a potentially different rate, giving rise to a competition of chase and escape. 

We show well-definedness of the model, isolate regimes of global survival and extinction of the infection and present estimates on local survival. The proofs rest on comparisons to the process on trees, percolation arguments and finite-degree approximations of the underlying random graphs. 
\end{abstract}

\section{Setting and main results}
\begin{figure}[!htpb]
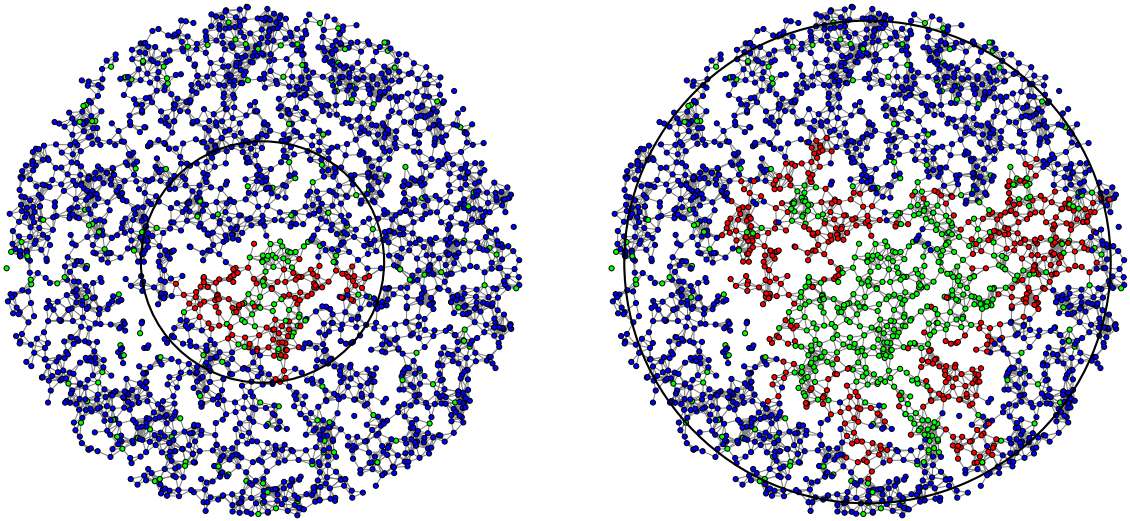

\centering
\begin{subfigure}{.5\textwidth}
  \centering
\include{Pix-WK-Stop-Exp-R_1}
\end{subfigure}%
\begin{subfigure}{.5\textwidth}
  \centering
\include{Pix-WK-Stop-Exp-R_2}
\end{subfigure}
\caption{Realization of a random network of nodes that are either infected (red), susceptible (blue) or white knights (green), in a disc at two finite times (left and right). Edges (gray) connect nodes at close proximity and allow for transmission of malware or patching. Circles (black) indicate the maximal distance of the malware to the origin in which the malware was placed initially.}
\label{Pix-WK-Stop-Exp}
\end{figure}

In this paper, we pick up a line of research, that very recently has attracted some attention, about the survival of some species when chased by another species, see~\cite{DuJuTa18} and references therein. To add another interpretation, our motivation for the model stems from applications in device-to-device networks. Imagine a device is infected by some malware at time zero, where the device is a vertex in some random geometric graph representing an ad-hoc telecommunication network. In order to stop the malware from spreading into the system like an infection, special devices can be introduced that have the ability to remove the malware from infected neighboring devices. The special devices that carry the patch are sometimes called \emph{white knights}. The white knights are not allowed to simply transfer the patch to any device in their vicinity, but only to malware-carrying devices. This is motivated by the fact that legal regulations do not allow forceful installation of patches without the consent of susceptible devices, unless the device poses a detected threat. However, once the safety hazard is detected, the operator is allowed to take countermeasures. Once the patch is installed, the infected device becomes a white knight itself, creating again a chase-escape dynamics where the malware is followed by white knights, see Figure~\ref{Pix-WK-Stop-Exp} for an illustration. We present more background in Section~\ref{sec-Stra}. 

\medskip
More specifically, we consider a random network of nodes in $\R^d$ given by a {\em homogeneous Poisson point process} $X=\{X_i\}_{i\in N}$ with intensity $\mu>0$, plus an additional node $o$ at the origin. Any two nodes $X_i,X_j\in X\cup\{o\}$ are connected by an edge if and only if $X_j\in B_r(X_i)$, where $B_r(x)$ denotes the ball centered at $x\in\R^d$ with radius $r>0$ that we treat as a fixed system parameter. This gives rise to the classical Boolean model or {\em Gilbert graph} $g_r(X\cup\{o\})$ from stochastic geometry. Any particle $X_i$ at time $t\ge 0$ can be in one of three states, 
\begin{align*}
\xi(t,X_i)=
\begin{cases} {\rm S}, &\text{if } X_i \text{ is {\it susceptible} at time }t, \\ 
{\rm I}, &\text{if }  X_i \text{ is {\it infected} at time } t,\\ 
{\rm W}, &\text{if } X_i \text{ is a {\it white knight} at time } t.
\end{cases}
\end{align*} 
For the set of all susceptible, infected and white-knight nodes at time $t$, we write respectively 
\begin{align*}
S(t)&= \{X_i \in X\cup\{o\} \colon \xi(t,X_i)={\rm S}\},\\
I(t)&= \{X_i \in X\cup\{o\} :\,  \xi(t,X_i)={\rm I}\} \text{ and} \\
W(t)&=\{X_i \in X\cup\{o\} :\,  \xi(t,X_i)={\rm W}\}.
\end{align*}
The propagation mechanism is given by a continuous-time Markov jump process with the following transition rates. 
\begin{enumerate}
\item If $\xi(t,X_i)={\rm S}$, then $X_i$ becomes infected with rate $\laI\#\big(I(t)\cap B_r(X_i)\big)$
 and
 \item if $\xi(t,X_i)={\rm I}$, then $X_i$ becomes a white knight with rate $\laW\#\big(W(t)\cap B_r(X_i)\big)$,
\end{enumerate}
where $\laW>0$ is the \emph{patch rate} and $\laI> 0$ is the \emph{infection rate}.  Through a re-scaling on the time axis we can set without loss of generality~$\laW=1$ for the remainder of this manuscript. 

\medskip
We will always assume that the infection starts at the origin, i.e., $I(0)=\{o\}$. Note here that putting an additional node into the network at the origin, amounts to considering the network under the Palm distribution, by Slivnyak--Mecke's theorem. Using this interpretation, the initially infected node is a 'typical node' in the system. As for the initial configuration of white knights, we assume them to be an i.i.d.~thinning of $X$ with parameter $p\in[0,1]$. In particular, the initial set of white knights $\XW=W(0)$ then is again a Poisson point process with intensity $\muW=p\mu$ and the initial set of susceptible particles $\XS=S(0)$ is a Poisson point process with intensity $\muS=(1-p)\mu$. It is one of the classical result in continuum percolation theory that there exists a unique critical intensity $0<\muR^{\rm c}<\infty$ such that for $\muS>\muR^{\rm c}$, the graph $g_r(\XS \cup \{o\})$ contains a unique infinite component of nodes with probability one and for $\muS<\muR^{\rm c}$, the graph $g_r(\XS \cup \{o\})$ contains no infinite component of nodes with probability one.

\medskip
More formally, let us denote by $\P=\Pspa \otimes \Pinf$ the joint probability distribution of the network model $\Pspa$ and the propagation model $\Pinf$ with initial configuration $\xi$. More precisely, by $\Pspa$ we denote the joint distribution of the superposition of the independent Poisson point processes $\XS$ and $\XW$, with an additional node at the origin. 
For a given realization $g_r(\XS\cup \XW\cup\{o\})$ of the associated Gilbert graph, $\Pinf$ denotes the probability kernel of the Markov propagation model with $\xi$ the initial configuration of states of nodes. $\xi$ is given by 
\begin{align*}
\xi(X_i) =\xi(0,X_i)= 
  \begin{cases}
    \rm{I},  & \text{for } X_i = o, \\
    \rm{S}, & \text{for } X_i \in \XS, \\
   \rm{W}, & \text{for } X_i \in \XW.
  \end{cases}
\end{align*}
Our first result establishes well-definedness of the propagation model for almost-all network realizations. The proof is presented in Section~\ref{sec-Proofs}.
\begin{prop}[Well-definedness]\label{prop-well}
$\Pinf$ is a well-defined standard continuous-time Markov jump process on $\{{\rm I,S,G}\}^{\XS\cup\XW\cup\{o\}}$, $\Pspa$-almost surely. 
\end{prop}
Our main interest lies in the analysis of extinction and survival of the infection. 
We denote by 
\begin{align*}
\EExt&= \{\text{there exists } t\ge 0 \colon  \#I(t)=0 \}, 
\end{align*}
the event of \emph{extinction} of the infection and call $\ESurv$ the \emph{survival} event. 
Due to the percolation properties of the underlying network, we further want to distinguish two types of survival. Let us call 
\begin{align*}
\ELoc &= \{\text{for all } t\ge 0 \colon \#I(t)>0\}\cap\{\#\Big(\bigcup_{t\ge 0}I(t)\Big)<\infty\},
\end{align*}
the event of \emph{local survival}, in which infected nodes are present for all times, but the number of such nodes is finite. On the other hand, let 
\begin{align*}
\EGlob&= \{\text{for all } t\ge 0 \colon \#I(t)>0\}\cap\{\#\Big(\bigcup_{t\ge 0}I(t)\Big)=\infty\},
\end{align*}
denote the event of \emph{global survival}, in which the infection never disappears and additionally reaches infinitely many nodes. All the events $E$, $L$ and $G$ of course depend on all the model parameters, which we suppress in the notation for convenience. 

\medskip
Let $\k_r$ denote the Lebesgue volume of the ball $B_r(o)$ and assume $\muS\k_r\ge 1$, then we define the quantity
$$\rho(\muS\k_r) = 2\muS\k_r-1-2\sqrt{(\muS\k_r)^2-\muS\k_r},$$
and note that $\rho:\, [1,\infty)\to (0,1]$, $x\mapsto\rho(x)$ is strictly decreasing. Further note that if $\muS\ge \muR^{\rm c}$, then $\muS\k_r\ge 1$, see for example~\cite[Equation 6.2]{pcPerc}. 
We are now in the position to state our main result about extinction. 

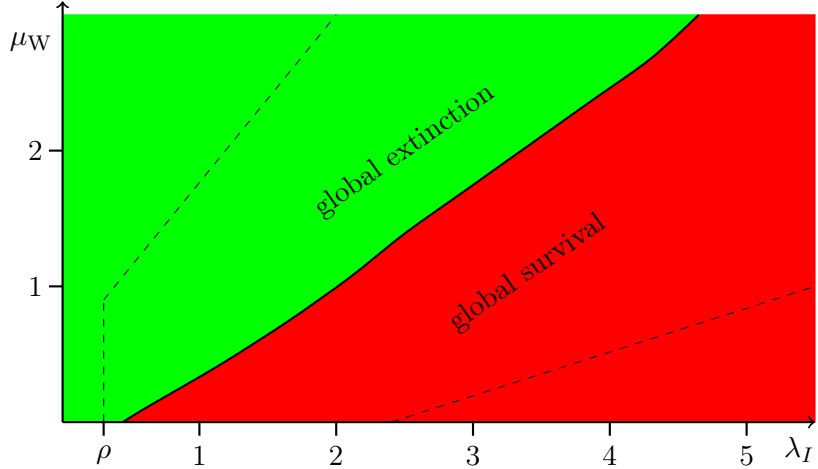
\begin{figure}[!htpb]
\centering
 \begin{tikzpicture} [xscale=1.8,yscale=1.8]

 \draw [fill, red] (0.44,0) --  (0.6,0.1) -- (1.1500000000000004,0.42222222222222217) -- (1.6500000000000008,0.7444444444444444) -- (2.1000000000000005,1.0666666666666667)-- (2.4999999999999996,1.3888888888888888) -- (2.9499999999999984,1.711111111111111) -- (3.3999999999999972,2.033333333333333) -- (3.849999999999996,2.3555555555555556) -- (4.299999999999995,2.6777777777777776) -- (4.649999999999994,3.0) -- (5.5,3.0) -- (5.5,0.0) -- cycle;

\draw [fill, green] (0.44,0) -- (0.6,0.1) -- (1.1500000000000004,0.42222222222222217) -- (1.6500000000000008,0.7444444444444444) -- (2.1000000000000005,1.0666666666666667)-- (2.4999999999999996,1.3888888888888888) -- (2.9499999999999984,1.711111111111111) -- (3.3999999999999972,2.033333333333333) -- (3.849999999999996,2.3555555555555556) -- (4.299999999999995,2.6777777777777776) -- (4.649999999999994,3.0) -- (0,3.0) -- (0.0,0.0) -- cycle;

  \draw [thick, black, xshift=0cm]  
		 plot [smooth] coordinates 	
{(0.44,0)(0.6,0.1)(1.1500000000000004,0.42222222222222217)(1.6500000000000008,0.7444444444444444)(2.1000000000000005,1.0666666666666667)(2.4999999999999996,1.3888888888888888)(2.9499999999999984,1.711111111111111)(3.3999999999999972,2.033333333333333)(3.849999999999996,2.3555555555555556)(4.299999999999995,2.6777777777777776)(4.649999999999994,3.0)};
 
		\draw[thick, ->] (0,0) -- (5.5,0);
		\draw[thick, ->] (0,0) -- (0,3.1);
		\coordinate[label = -180 : $\muW$](a) at  (0.0,2.8);
		\coordinate[label = 0 : $\lambda_I$ ](a) at  (5.2,-0.2);
\foreach \x in {1,2,3,4,5}
 \node[anchor=north] at (\x,-0.1) {\x} ;
 \node[anchor=north] at (0.3,-0.1) {$\rho$} ;

\foreach \x in {0.3,1,2,3,4,5}
\draw[thick] (\x,0) -- (\x,-0.1);

\foreach \y in  {1,2}
 \node[anchor=east] at (-0.07,\y) {\y};

\foreach \y in  {1,2}
 \draw[thick] (-0.10,\y) -- (0,\y);

\draw[black, dashed] (0.3,0.9) -- (2,3);
\draw[black, dashed] (0.3,0.0) -- (0.3,0.9);

\draw[black, dashed] (2.4,0.0) -- (5.5,1);

\node[rotate=35] at (2.5,2.0) {global extinction};
\node[rotate=35] at (3.4,1.1) {global survival};

\end{tikzpicture}
\caption{Simulated phase diagram of global survival and extinction in the plane of infection rate vs.~white-knight intensity for $d=2, r=1,\muS=3,\laW=1$. Regimes where extinction and survival are covered by Theorem~\ref{thm-GloExt} and Theorem~\ref{thm-GloSur} are roughly indicated by dashed lines. 
}
\label{Pix-PhaseDiagram}
\end{figure}

\begin{thm}[Global extinction]\label{thm-GloExt}
If $0\le\muS<\muR^{\rm c}$, then $\P(G)=0$ for all $\laI\ge 0$ and $\muW\ge0$. Further, if $\muS\ge \muR^{\rm c}$ and $\laI\ge 0$, then there exists $\muW^{\rm c}(\laI,\muS)<\infty$ such that for all $\muW>\muW^{\rm c}(\laI,\muS)$ we have that $\P(G)=0$. Finally, if $\muS\ge \muR^{\rm c}$ and $\laI\le\rho(\muS\k_r)$, then $\muW^{\rm c}(\laI,\muS)=0$.
\end{thm}
Before we present our result about global survival, let us comment on the preceding theorem. In simple terms, Theorem~\ref{thm-GloExt} says that if the graph of susceptible nodes is insufficiently connected, i.e., $\muS<\mu_r^{\rm c}$, then global survival is impossible for any infection rate and even without any white knights in the system. Next, if the infection is too weak with respect to the intensity of susceptible nodes, i.e., $\laI\le \rho(\muS\k_r)$, then global survival is also impossible for any positive intensity of white knights. On the other hand, for any infection rate, sufficiently many white knights in the system lead to global extinction, see the green part in Figure~\ref{Pix-PhaseDiagram} for an illustration. 
As we will explain later in Section~\ref{sec-Stra}, the specific form of the threshold $\rho$ is due to a comparison to the propagation model on trees. We also present a more detailed explanation and references to preceding research there.

\medskip
The next result is about global survival, see the red part in Figure~\ref{Pix-PhaseDiagram} for an illustration.

\begin{thm}[Global survival]\label{thm-GloSur}
For all  $\muS>\muR^{\rm c}$ and $\muW\ge 0$ there exists $\lambda_{\rm I}^{\rm c}(\muW,\muS)<\infty$ such that for all $\lambda_{\rm I} >\lambda_{\rm I}^{\rm c}(\muW,\muS)$, we have that $\P(G) > 0$.
\end{thm}
In words, Theorem~\ref{thm-GloSur} states that there is a positive chance for global survival if the underlying graph of susceptible nodes is sufficiently connected and the infection is strong enough to overcome the chasing white knights. 

\medskip
For statements about local survival, let us introduce the notation $C_o$ for the cluster of all nodes in $X=\XS\cup \XW$ for which there exists a path in  $g_r(X\cup \{o\})$ connecting them to the origin. Similarly, we denote by $C_o^S$ the cluster of nodes connected to the origin in $g_r(\XS\cup \{ o \})$. Further, denote by $\theta(\muS)=\Pspa(\# C_o^S= \infty)$, the percolation probability of the process of susceptible nodes. Note that the infection can never leave the set $C_o^S$ and hence global survival is impossible if $C_o^S$ is finite and this implies $\P(G)\le\theta(\muS)$. 
\begin{lem}[Local survival]\label{lem-LocSur}
For all parameters,
$$\exp(-\mu\k_r)\le \P(L)\le (1-\theta(\muS))\exp(-\muW\k_r).$$  
\end{lem}
In particular, if the process of susceptible nodes is subcritical, we have $\P(G)=0$, $\P(L)\le\exp(-\muW\k_r)$ and thus $\exp(-\mu\k_r)\le \P(E^c)\le\exp(-\muW\k_r)$. If the process of susceptible nodes is supercritical but the other parameters guarantee global extinction as described in Theorem~\ref{thm-GloExt}, then for the survival probability $\exp(-\mu\k_r)\le \P(E^c)\le(1-\theta(\muS))\exp(-\muW\k_r)$. Note that $\mu\mapsto\theta(\mu)$ tends to one exponentially fast, see~\cite{uniqFin}, and hence local survival is exponentially unlikely for dense networks both in $\muS$ and $\muW$.

\medskip
In the next section we explain the strategy of the proofs, and comment on related results in the literature. The proofs are presented in Section~\ref{sec-Proofs}. 


\subsection{Acknowledgement}
We thank the team from Orange S.A., in particular Ali--Malek Boubaya as well as Wolfgang K\"onig and Andr\'as T\'obi\'as for inspiring discussions. This work was funded by the German Research Foundation under Germany's Excellence Strategy MATH+: The Berlin Mathematics Research Center, EXC-2046/1 project ID: 390685689 as well as Orange Labs S.A..

\section{Strategy of proofs}\label{sec-Stra}
The study of epidemic models defined in terms of interacting particle systems with some additional randomness coming from an environment has, by now, a long history, see for example the early works~\cite{Li92,An92}. The consideration of such processes on random graphs has attracted attention more recently, see for example~\cite{Du10a,Du10b} and references therein. Here, a particularly interesting class is given by random graphs with a prescribed degree distribution, see for instance~\cite{ChDu09,MoVaYa13,MoMoVaYa16,MoVa16}. The literature on interacting particle systems on random geometries, i.e., where the random graphs are embedded in space, and in particular do not obey any degree bounds, is much sparser and younger. Notable here is the work~\cite{MeSi16}, which establishes existence of a subcritical phase of the contact process on Gilbert graphs and Poisson--Delaunay tessellations. 

\medskip
As mentioned in the beginning, propagation models analogue to the one presented in this manuscript have been studied by Lalley, Tang, Kordzakhia as well as Durrett and coauthors in recent years on various fixed networks such as trees~\cite{Ko05} and lattices~\cite{DuJuTa18} also via simulations, see~\cite{TaKoLa18}, and mainly motivated by applications in probabilistic biology.  
Apart from our generalization towards random geometries, another difference in the analysis presented here is that, in the initial configuration, the white knights form a Poisson point process of infinitely many nodes, whereas in the preceding works on fixed geometries there is only a finite number of white knights present in the system. 
Nevertheless, our proofs are partially based on the results for fixed networks, in particular the tree considered by Kordzakhia~\cite{Ko05}, since there it is possible to derive explicit bounds for the infection rate by balancing numbers of paths compared with propagation along one path, which is then in fact one dimensional. More precisely, let us consider our propagation model on a fixed connected graph $H$ that includes a root $\{o\}$ and a generic point $\{o'\}$ which is only connected to the root by a single edge. Let the starting configuration be given by 
\begin{align*}
\xi'(x) =
  \begin{cases}
    \rm{W},  & \text{for } x  = o',  \\
    \rm{I},  & \text{for } x  = o,  \\
   \rm{S}, & H\setminus \{o,o'\}.
  \end{cases}
\end{align*}
Let $\laI^{\rm c}(H)$ denote the critical rate for (global) extinction of the infection based on the propagation model as explained above, where the underlying Gilbert graph and the initial condition $\xi$ are replaced by $H$ and $\xi'$. Then, the following result is proved in~\cite[Theorem 1 and Corollary  2]{DuJuTa18}.
\begin{lem}[Extinction on fixed networks] \label{Cor2Durrett}
Let $\Gamma_n(H)$ denote the set of self-avoiding paths of length $n$ in $H$, starting from the root. If there exists $k\in\{2,3,\dots\}$ such that for all $n\in\N$ we have that $\#\Gamma_n(H) \le k^n$, then 
$$\laI^{\rm c}(H) \ge\laI^{\rm c}(\mathbbm{T}_k)= 2k-1-2\sqrt{k^2-k},$$ 
where $\mathbbm{T}_k$ is the rooted $k$-ary tree.
\end{lem}
The next result establishes existence of a constant, the {\em connective constant}, bounding the number of self-avoiding paths for almost-all realizations of the Gilbert graph. 
\begin{lem}[Connective constant]\label{lem-ConnectiveConstant}
For all $\muS, r >0$ and all $\gamma>\muS\k_r$ we have 
\begin{align}
\limsup_{n\uparrow\infty}n^{-1}\log\#\Gamma_n\big(g_r(\XS\cup\{o\})\big)\le \log \gamma,
\end{align}
for $\Pspa$-almost all $\XS$. 
\end{lem}
We present the proof of Lemma~\ref{lem-ConnectiveConstant} in Section~\ref{sec-Proofs}. Now the last statement of Theorem~\ref{thm-GloExt} is an immediate consequence of the following proposition,  which leverages Lemma~\ref{lem-ConnectiveConstant} and the proof idea of Lemma~\ref{Cor2Durrett} to the setting of infinitely-many white knights. 
\begin{prop}\label{prop-GloExt1} 
If $\muS\ge \mu_r^{\rm c}$, $\laI\le \rho(\muS\k_r)$ and $\muW>0$, then $\P(G\cap \{\# C_o=\infty\})=0$. 
\end{prop}

\begin{proof}[Proof of Theorem~\ref{thm-GloExt}, first and last statement]
First of  all, if $\muS<\mu_r^{\rm c}$, then, as mentioned above, $\P(G)\le\theta(\muS)=0$. For the last statement, let $\muS\ge\mu_r^{\rm c}$, $\laI\le \rho(\muS\k_r)$ and $\muW>0$, then using Proposition~\ref{prop-GloExt1}, we have that 
\begin{align*}
\P(G)=\P(G\cap \{\# C_o=\infty\})+\P(G\cap \{\# C_o<\infty\})=0,
\end{align*}
since global survival is impossible on finite clusters. 
\end{proof}
The technique used for the proof of Proposition~\ref{prop-GloExt1}, which is based on the works of Durrett and coauthors, fail in the regimes of large $\laI$, since the discovery of a new white knight on a one-dimensional path not only stops the infection on this path, but also creates new white knights. However, with the help of percolation arguments we prove in Section~\ref{sec-Proofs} the following proposition, which immediately implies the second statement of Theorem~\ref{thm-GloExt}. 
\begin{prop}\label{prop-GloExt2}
For all $\laI\ge 0$ and $\muS\ge0$, there exists $\muW^{\rm c}(\laI,\muS)<\infty$ such that for all $\muW>\muW^{\rm c}(\laI,\muS)$, we have that $\P(G\cap\{ \#C_o= \infty\})=0$.
\end{prop}
\begin{proof}[Proof of Theorem~\ref{thm-GloExt}, second statement] 
In the setting of Proposition~\ref{prop-GloExt2}, we have $\P(G)=\P(G\cap \{\# C_o=\infty\})+\P(G\cap \{\# C_o<\infty\})=0$.
\end{proof}
Before we comment on the proof of Theorem~\ref{thm-GloSur}, let us give a heuristic for the linear dependence of $\laI\mapsto\muW^{\rm c}(\muS,\laI)$ indicated in Figure~\ref{Pix-PhaseDiagram}.
For this consider the rescaled process of white knights with intensity $\alpha \muW$ and patch rate $1/\a$ for $\a>0$. Then, the rescaled chase-escape model, in the limit as $\alpha$ tends to infinity, converges to a {\em contact process} on $g_r(\XS\cup\{o\})$ with recovery rate given by $\muW\k_r$ and infection rate $\laI$. 
Then, the linear dependence in $\laI$ emerges through a re-scaling on the time axis.

\medskip
The main challenge for the proof of Theorem~\ref{thm-GloSur} lies in the fact that almost-all network realizations have an unbounded degree. Although large degrees should support survival of the infection, lack of monotonicity prevents us from using this idea directly. See also our comments on monotonicity below. However, we can estimate the network by graphs of bounded degree and use discrete percolation arguments for the approximations. The proof of Theorem~\ref{thm-GloSur} is presented in Section~\ref{sec-Proofs}. 

\medskip
Finally the quantitative statements about local survival as presented in Lemma~\ref{lem-LocSur}, are a simple consequence of void-space probabilities, let us give the proof here as well.
\begin{proof}[Proof of Lemma~\ref{lem-LocSur}] 
First, note that local survival is only possible if $\# C^S_o<\infty$ since otherwise, with probability one, there exists a white knight in the set $B_r(C^S_o)=\bigcup_{X_i\in C^S_o}B_r(X_i)$, which is eventually reached by the infection. Consequently, the infection can only survive if it escapes towards infinity on $C^S_o$. In particular, local survival is only possible if no white knights are in $B_r(C^S_o)$, i.e., 
\begin{align*}
\P(L)&=\P(L\cap\{\# C^S_o<\infty\})=\E\big[\exp(-\muW|B_r(C^S_o)|)\one\{\# C^S_o<\infty\}\big]\\
&\le (1-\theta(\muS))\exp(-\muW\k_r). 
\end{align*}
On the other hand, the infection certainly survives if the origin is isolated. This gives the lower bound. 
\end{proof}

\begin{figure}[!htpb]
\centering
\begin{tikzpicture}[scale = 2.0]
		\draw[thick] (1,0) -- (3.5,0);
		\draw[thick, dashed] (-1,0) -- (1,0);
		\draw[fill, green] (-1,0) circle (7pt);
		\draw[fill, white] (0,0) circle (7pt);		
		\draw[dashed](0,0) circle (7pt);
		\draw[fill, red] (1,0) circle (7pt);
		\draw[fill, blue] (2,0) circle (7pt);
		\draw[fill, blue] (3,0) circle (7pt);
		\draw[fill, black] (3.5,0) circle (1.5pt);
		\draw[fill, black] (3.7,0) circle (1.5pt);
		\draw[fill, black] (3.9,0) circle (1.5pt);
\end{tikzpicture}
\caption{Illustration of a configuration of a white knight (green, position $1$), an infected node (red, position $3$) and susceptible nodes (blue) for which adding an infected node at position $2$ would stop the infection that would otherwise propagate to infinity.}
\label{Pix_non_monotone}
\end{figure}
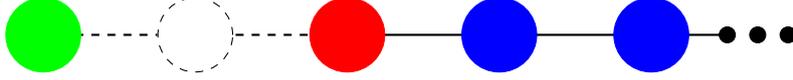

Let us finish this section by commenting on monotonicity properties of the phase diagram as sketched in Figure~\ref{Pix-PhaseDiagram}. Both, simulations and common sense suggest existence of a unique phase-separating curve and several monotonicities depending on the parameters. For example that additional infected nodes or an increase in the infection rate should increase the probability for the infection to survive. 
However, to prove existence, uniqueness and monotonicities is challenging, mainly because of the existence of configurations that exhibit counterintuitive effects, standing in the way of coupling arguments. 
Let us give one example here. Note that white knights can only act towards perviously infected nodes, and therefore an increase in infected nodes also benefits the spread of white knights. To illustrate this, consider the nearest-neighbor graph on $\N$ as presented in Figure~\ref{Pix_non_monotone} with a white knight at node $1$, an infection at node $3$ and all other nodes being susceptible. Imagine node $2$ and its associated edges were absent, then the infection would spread towards infinity unstopped for any positive infection rate. Now, if we add an infected node at position $2$, then more infections are in the system. Still, if $\laI<1$, the infection will now go extinct. Other examples can be constructed to also showcase configurations where an increase of $\laI$ leads to a decrease for the probability of survival of the infection. 

\medskip
Let us finally mention that in the survival regime it is reasonable to believe that there exists an infection speed $\a$ depending on all the parameters of the system, such that as $t \rightarrow \infty$ we have $B_{\alpha t(1-\epsilon)}(o)\subset I(t) \subset B_{\alpha t(1+\epsilon)}(o)$ with probability one, conditioned on the event that the origin is connected to infinity. 

\medskip
In the following section we present all remaining proofs. 

\section{Proofs}\label{sec-Proofs}
In  the sequel, we denote by $\Espa$ and $\Einf$ the expectations associated to $\Pspa$ and $\Pinf$, respectively.  
We abbreviate for balls, $B_n(o)$ by $B_n$ and for boxes, $Q_n(o)$ by $Q_n$. For any $A\subset\R^d$ we write $A^{\rm c}=\R^d\setminus A$. 

\begin{proof}[Proof of Proposition~\ref{prop-well}]
Note that for positive intensities, $\Pspa$-almost surely, all the Gilbert graphs $g_r(\XS\cup\{o\})$, $g_r(\XW\cup\{o\})$ and $g_r(X\cup\{o\})$ have an unbounded degree, potentially leading to blow-ups in finite time for the propagation model.  In particular, well-definedness of the propagation model on the Gilbert graphs can not be guaranteed using the standard conditions based on generators. For example there exists no finite bound on the transition rate uniformly in the nodes, as required in~\cite[Proposition 3.2 Chapter 3]{Li85}.

However, in our case, we can establish well-definedness due to the fact that in our initial condition there is only one infection present. More precisely, consider the process $(\xi_n(t,\cdot))_{t\ge 0}$ on $g_r\big((X\cup\{o\})\cap B_n\big)$ with $n\in r\N $, defined via the same rates as presented above. 
Then, for $\Pspa$-almost-all network realizations, the process $(\xi_n(t,\cdot))_{t\ge 0}$ is well-defined and standard as a finite state space Markov process. Next, let 
$$\tau_n= \inf\{t>0\colon  I_{n+1}(t) \not\subset B_n\}$$
denote the time at which the infection, based on the process $\xi_{n+1}(t,\cdot)$, hits the boundary of $B_n$. 
Then, for all $t< \tau_n$ our original process $\xi(t,\cdot)$ coincides with $\xi_n(t,\cdot)$ and is thus well-defined.
 
What remains to be shown is that $\P$-almost surely $\lim_{n \uparrow \infty} \tau_n = \infty$.
To show this, we will use a large-deviation argument to establish a bound on the minimal asymptotic speed of the infection process and apply the Borel--Cantelli lemma. More precisely, first note that for any $\a>0$ and $n\in r\N$, we have that
\begin{align}\label{upper-bound_1}
\Pinf \left(\tau_n < n/\alpha \right) 
\le \#\Gamma_{n/r}\Pinf\left(S_{n/r}<n/\alpha \right), 
\end{align}
where $\#\G_m$ is the number of self-avoiding paths in $g_r(\XS\cup\{o\})$ of length $m\in \N$ starting in $\{o\}$ and $S_m$ is the sum of $m$ independent exponentially-distributed random variables with parameter $\laI$. For the estimate~\eqref{upper-bound_1} we used that, in order for the infection to reach $B_n^{\rm c}$ in less than $n/\a$ time, at least $n/r$ infection events happened along at least one of the self-avoiding paths of length $n/r$ in less than $n/\a$ time. Now, using Lemma~\ref{lem-ConnectiveConstant} and large-deviation bounds for independent exponentially-distributed random variables, we obtain for $\alpha > r\laI$ and $\Pspa$-almost all network realizations the bound 
\begin{align}\label{eq_prop1.1_2}
\limsup_{m \uparrow \infty} \frac{1}{m} \log \left(\#\Gamma_{m}\Pinf \left( S_{m}<\frac{mr}{\alpha} \right)\right)
 	\le \log\gamma -\frac{\laI r}{\alpha } +1+ \log\frac{\laI r}{\alpha }= -C_{\g,\laI,r}(\alpha).
\end{align}
In particular, for $\a>\a_{\rm c}$ with $\alpha_{\rm c}= \inf\{\alpha >r\laI \colon -C_{\g,\laI,r}(\alpha)>0\}$,
there exists $n_o\in r\N$ such that 
\begin{align*}
\sum_{n \in r\N} \Pinf \left(\tau_n < n/\alpha \right)  < n_o+\sum_{n\in r\N\colon n> n_o}\exp(-nC_{\g,\laI,r}(\alpha)/2) < \infty.
\end{align*}
Finally, by an application of the Borel--Cantelli lemma, $\Pspa$-almost surely, we have that \linebreak $\Pinf \left(\limsup_{n\uparrow\infty}\{\tau_n < n/\alpha\} \right)=0$ and thus almost surely under $\P$, for all but finitely many $n$, we have $\tau_n \ge n/\alpha$, which finishes the proof. 
\end{proof}
Let us note that the above proof of existence of the process also derives bounds on the minimal speed of propagation of the infection in terms of solutions to fixed-point equations determining the critical speed $\a_{\rm c}$.

\begin{proof}[Proof of Lemma~\ref{lem-ConnectiveConstant}]
First, we compute the expectation of $\#\Gamma_n$ by multiple applications of Slivnyak--Mecke's theorem. More precisely, let $X_{i_0}$ be an alternative notation for $o$, then for any $n\in\N$ we have that 
\begin{align*}
\Espa[\#\Gamma_n]&=\Espa\Big[\sum_{i_1,\dots,i_n\in N}^{\neq}\prod_{k=1}^n \one\{|X_{i_k}-X_{i_{k-1}}|<r\}\Big]\\
&=\muS\int\Espa\Big[\sum_{i_1,\dots,i_{n-1}\in N}^{\neq}\prod_{k=1}^{n-1} \one\{|X_{i_k}-X_{i_{k-1}}|<r\}
\one\{|x-X_{i_{n-1}}|<r)\}\Big]\,\d x\\
&=\muS\k_r\,\Espa\Big[\sum_{i_1,\dots,i_{n-1}\in N}^{\neq}\prod_{k=1}^{n-1} \one\{|X_{i_k}-X_{i_{k-1}}|<r\}
\Big]\\
&=(\muS\k_r)^{n},
\end{align*}
where the $\neq$ sign indicates that we sum over mutually distinct indices. Next, using the Markov inequality, we have
\begin{align*}
\Pspa(\#\Gamma_n\ge \gamma^n) \le \g^{-n}\Espa[\#\Gamma_n]= \big(\muS \k_r/\gamma\big)^n,
\end{align*}
and thus, for $ \gamma > \muS\k_r$ we obtain  $\sum_{n \in \N}\Pspa(\#\Gamma_n\ge \gamma^n) < \infty$. Hence, the Borel--Cantelli lemma yields that $\#\Gamma_n> \gamma^n$ only for finitely-many $n$. This completes the proof. 
\end{proof}

\begin{proof}[Proof of Proposition~\ref{prop-GloExt1}]
It suffices to consider the critical and supercritical regimes where $\muS\ge\muR^{\rm c}$.
Let us abbreviate $A=\bigcup_{t\ge 0} I(t)$. First, note that it suffices to prove that for almost-all network realizations we have that $E^\xi[\# A]<\infty$ for $\laI\le\r(\muS\k_r)$ and any $\muW>0$, as this implies $\Pinf(G)=0$. For this, we want to bound $\#A$ from above by numbers of infected nodes along one-dimensional paths. We have to be careful here since by the lack of monotonicity, as explained at the end of Section~\ref{sec-Stra}, we cannot simply remove white knights in order to produce a situation as in Lemma~\ref{Cor2Durrett}. However, note that for every $x\in A$ there must exist a self-avoiding path $\phi_n(x)\subset g_r(\XS\cup\{o\})$ of some finite length $n=n(x)$, started at the origin, such that $x$ received its infection along $\phi_n(x)$. Let $I_n\subset A$ denote the set of nodes in $A$ that were infected at some time along a path of length $n$. Then, for $\Pspa$-almost all $X$ with $\XW\cap B_r(o)\neq\emptyset$, i.e., realizations of the network where the origin is directly adjacent to a white knight, we have
\begin{equation}\label{DurInc}
\begin{split}
E^\xi[\#A]&\le \sum_{n\ge 0}E^\xi[\# I_n]\le \sum_{n\ge 0}\#\G_n P'(o\text{ infects }n \text{ on }\N\cup\{o\}\cup\{o'\}),
\end{split}
\end{equation}
where $P'$ denotes the distribution of the process described in Lemma~\ref{Cor2Durrett} for $H=\N\cup\{o\}\cup\{o'\}$. 
Here, we used that more neighboring white knights along a fixed path lead to even smaller probability of survival.

An application of Lemma~\ref{Cor2Durrett}, as presented in~\cite[Theorem 1 and Corollary  2]{DuJuTa18}, yields that 
\begin{align*}
P'(o\text{ infects }n \text{ on }\N\cup\{o\}\cup\{o'\})\le C(\laI)\laI'^n n^{-3/2},
\end{align*}
for $\laI'=4\laI/(1+\laI)^2$, which is derived via the reflection principle in case of $\laI<1$. Note that by assumption, $\g>1$ and thus for $\laI\le  2\gamma - 1 - 2\sqrt{\gamma^2-\gamma}<1$ in particular $\laI'<1$. Hence, for such $\laI$, the right-hand side of~\eqref{DurInc} is finite and thus $E^\xi[\#A]<\infty$ for $\Pspa$-almost all $X$ with $\XW\cap B_r(o)\neq\emptyset$.

To finish the proof let us consider the network realizations where $\XW\cap B_r(o)=\emptyset$, i.e., where the origin is not adjacent to a white knight. Denote the set of all connected finite subsets of nodes that contain at most one node adjacent to a white knight by 
$$\mathcal{J}= \big\{ J\subset \XS \cup \{o\}\colon o \in J, \#J < \infty, J \text{ connected},\#\{X_i\in J\colon \dist(X_i,X_W)<r\}=1 \big\},$$
where $\dist(x,B)=\inf\{|x-y|\colon y\in B\}$ denotes the distance of a point $x\in \R^d$ to a set $B\subset\R^d$. For any $J\in\mathcal J$ we define $\xi^J$ to be the configuration where the infected nodes are precisely given by $J$, i.e., we define for all $X_i\in X$
\begin{align*}
\xi^J(X_i)=\begin{cases}
    \rm{I},  & \text{if } X_i \in J  \\
    \rm{S}, & \text{if } X_i \in \XS \setminus J \\
    \rm{W},  & \text{if } X_{i}  \in \XW.    
  \end{cases}
\end{align*}
Then, for $\Pspa$-almost-all $X$ with $\# C_o = \infty$ we define $\tau= \inf\{t\ge 0\colon \xi(t,\cdot)=\xi^J\text{ for some }J\in\mathcal J\}$, the stopping time at which the process explores the white knights for the first time. With these definitions, note that
\begin{align*}
\sum _{J \in \mathcal{J}} \Pinf \left(\xi(\tau) = \xi^J\right) = 1,
\end{align*}
since every realization of the infection propagation on an infinite cluster eventually reaches a white knight in finite time. Furthermore, the node that explored the white knight is unique and will be denoted by $o_J$.
Therefore, for $\Pspa$-almost-all $X$ with $\# C_o = \infty$ we have
\begin{align*}
\Pinf(G) = \sum_{J \in \mathcal{J}}\Pinf \left(G \mid \xi(\tau)= \xi^J\right) \Pinf( \xi(\tau)= \xi^J).
\end{align*}
Now, due to the strong Markov property, we have 
\begin{align*}
\Pinf(G \mid \xi(\tau)= \xi^J)=P^{J} (G \mid \xi(\tau)=\xi^J)\le  P^{J} (G) /  P^{J} (\xi(\tau)= \xi^J) = 0,
\end{align*}
where $P^{J}$ denotes the infection process not started in $o,$ but in $o_J$, constructed on the same probability space as $\Pinf$.
The last equality holds, as $P^J( \xi(\tau)= \xi^J)>0$ and $P^J (G) = 0$, since with these definitions, there is a white knight next to the origin and the first part of the proof can be applied again.
Now, since $\mathcal{J}$ is countable, we can conclude that indeed for $\Pspa$-almost all $X$ we have $\Pinf(G)=0$. 
\end{proof}

\begin{proof}[Proof of Proposition~\ref{prop-GloExt2}]
The proof is based on a percolation argument. We call a node $X_i\in\XS$ an \emph{open node} if
\begin{enumerate}
\item $X_i$ is not isolated in $g_r(\XS\cup\{o\})$, and 
\item once $X_i$ is infected, then it transmits its infection towards at least one of its neighbors in $\XS$ (regardless if neighbors are already infected or patched) before $X_i$ is directly patched by a neighboring white knight in $\XW$.
\end{enumerate}
We call $X_i$ a \emph{closed node} otherwise. For example $X_i$ is an open node if it is not isolated in $g_r(\XS\cup\{o\})$ but isolated in $g_r(\XW\cup\{X_i\})$. Note that a node is labeled open or closed based on its neighborhood at initial time. Thanks to the strong Markov property, the probability of a node to be open does not depend on the time at which it becomes infected. Further note that an infinite self-avoiding path of open nodes in $g_r(\XS\cup\{o\})$ does not guarantee global survival of the infection since we do not require that the infection propagates to infinity. However, absence of an infinite self-avoiding path of open nodes implies absence of global survival in the realization. Indeed, if there is no path to infinity of nodes that are able to infect at least one neighboring node from the initially susceptible nodes, then in particular there is no path to infinity of nodes that are able to infect neighboring nodes that are still susceptible at the time at which the infection arrives. This comes from the fact that less susceptible nodes or more white knights in the neighborhood make it even harder to transmit the infection towards at least one neighboring susceptible node before being patched by a neighboring white knight. 

\medskip
In order to show absence of an infinite self-avoiding path of open nodes, let us discretize space into boxes $Q_{3r}(3rz)$ of side-length $3r$, centered at $3rz$ with $z\in\Z^d$ and use results on lattice percolation. We define
\begin{align*}
V_k=\{ x \in \R^d\colon \#(\XW\cap B_r(x) ) \ge k \},
\end{align*}
the set of all space points that have at least $k$ white knights in its neighborhood.
The site $z\in \Z^d$ is called an \emph{open site}, if one of the following events happens,
\begin{align*}
A_{m}(z)&=\{\#\big( \XS \cap Q_{3r}(3rz)\big)\ge m\},\text{ or}\\
B_{n}(z)&= \{\text{there exists } X_i \in \XS\cap  Q_{3r}(3rz)  \colon\#\big(\XS \cap B_r(X_i)\big)\ge n+1\},\text{ or}\\ 
C_{k,\muW}(z)&= \{Q_{3r}(3rz) \not\subset V_k\},\text{ or}\\
D(z)&=\{Q_{3r}(3rz)\text{ contains an open node}\}.
\end{align*}
Otherwise $z\in\Z^d$ is called a \emph{closed site}. Now, suppressing the dependence on $z$ if $z=o$, we have that 
\begin{align*}
\P(o\text{ is an open site})=1-\P(D^{\rm c}\mid A_m^{\rm c}\cap B_n^{\rm c}\cap C_{k,\muW}^{\rm c})\P(A_m^{\rm c}\cap B_n^{\rm c}\cap C_{k,\muW}^{\rm c}).
\end{align*}
If a node $X_i$ has $n\ge 1$ neighbors in $\XS$ and $k\ge0$ neighbors in $\XW$, the probability for $X_i$ to be a closed node is given by $k/(k+n\laI)$.
Hence, we can bound
\begin{align*}
\P(D^{\rm c}&\mid A_m^{\rm c}\cap B_n^{\rm c}\cap C_{k,\muW}^{\rm c})\\
&\ge \Espa\Big[\prod_{X_i\in \XS \cap Q_{3r}(o)}\frac{\#\big(\XW \cap B_r(X_i)\big)}{\#\big(\XW\cap B_r(X_i)\big)+\big(\#(\XS \cap B_r(X_i))-1\big)\laI}\mid A_m^{\rm c}\cap B_n^{\rm c}\cap C_{k,\muW}^{\rm c}\Big]\\
&\ge \Big(\frac{k}{k+n\laI}\Big)^{m},
\end{align*}
and 
\begin{align*}
\P(A_m^{\rm c}\cap B_n^{\rm c}\cap C_{k,\muW}^{\rm c})\ge 1-\big(\P(A_m)+\P(B_n)+\P(C_{k,\muW})\big).
\end{align*}
Now, for any $\laI\ge 0$ and $\eps>0$, there exist $m_o, n_o\in\N$ such that for all $m>m_o$ and $n>n_o$ we have 
\begin{align*}
\P(A_m)<\eps\qquad\text{ and }\qquad\P(B_n)<\eps
\end{align*}
and $k_o(m, n)\in\N$ such that for all $k>k_o(m,n)$ we have 
\begin{align*}
\Big(\frac{k}{k+n\laI}\Big)^{m}>1-\eps.
\end{align*}
Finally, we can then pick $\muW^{\rm c}(k_o(m,n))$ sufficiently large, such that for all $\muW>\muW^{\rm c}(k_o(m,n))$ also 
\begin{align*}
\P(C_{k,\muW})<\eps. 
\end{align*}
Together, the probability for an open site can be made arbitrarily small, since
\begin{align*}
\P(o\text{ is an open site})\le 1-(1-\eps)(1-3\eps)=4\eps-3\eps^2. 
\end{align*}

\medskip
The random field of good and bad sites constitutes a two-dependent site-percolation model on $\Z^d$ that can be dominated by an independent site-percolation model, using the domination-by-product-measures result~\cite[Theorem 0.0]{domProd}. Then, for sufficiently large $\muW$ we have absence of percolation of good sites. As the side-length of the boxes is larger than $2r$, this also excludes the possibility of continuum percolation of open nodes in the Gilbert graph, which then implies absence of global survival of the infection process also on the event $\{\# C_o=\infty\}$.
\end{proof}

\begin{proof}[Proof of Theorem~\ref{thm-GloSur}]
Similar to the proof of Proposition~\ref{prop-GloExt2} we will use a percolation argument for a suitably chosen discretization of $\R^d$ to show existence of an infinite cluster of nodes that transmit their infection faster than any cure attempt of neighboring white knights.  More precisely, we call $X_i\in\XS$ an \emph{open node} if $X_i$ transmits the infection to all its neighbors in $g_r(\XS)$ before an attempt to cure $X_i$ has been made by any neighbor in $g_r(\XS\cup\XW)$. We call $X_i$ a \emph{closed node} otherwise. With this definition, using also the strong Markov property, for the global survival of the infection, it suffices to prove existence of an infinite cluster of open nodes, which is connected to the origin, with positive probability.

We aim to achieve this by choosing the infection rate $\laI$ large, however, due to the unbounded degree of the graphs, there is no globally sufficiently large infection rate such that for $\Pspa$-almost all graphs the survival rate is above a fixed $\epsilon > 0$. We therefore introduce a discretization of $\R^d$ into boxes and distinguish those boxes in which the degrees are bounded. 
For this, let $g^{m,n}_{r,\laI}(\XS\cup \XW)$ denote the Gilbert graph with connectivity threshold $r>0$ and vertex set \begin{equation}\label{Survival_1}
\begin{split}
\{X_i\in\XS\colon &\#\big(\XW\cap B_r(X_i)\big)\le m,\, \#\big(\XS\cap B_r(X_i)\big)\le n+1\text{ and }X_i\text{ is an open node}\}.
\end{split}
\end{equation}
This is a Gilbert graph based on a dependently thinned Poisson point process, where vertices are removed if they have too many neighbors in $\XS$ and $\XW$ or transmit their infection too slowly. Let $\theta_{m,n,\laI}(\muS,\muW)$ denote the associated percolation probability. Our parameters are such that $\theta(\muS)>0$. We claim that 
\begin{align}\label{Continuity_1}
\lim_{m,n,\laI\uparrow\infty}\theta_{m,n,\laI}(\muS,\muW)=\theta(\muS).
\end{align}
If~\eqref{Continuity_1} holds, then this implies that for sufficiently large $m,n$ and $\laI$ also $g^{m,n}_{r,\laI}(\XS\cup\XW)$ is supercritical and the infection survives globally. 
In order to prove~\eqref{Continuity_1}, first note that for all $m,n,\laI$ we have that $\theta_{m,n,\laI}(\muS,\muW)\le \theta(\muS)$ since $g^{m,n}_{r,\laI}(\XS\cup \XW)$ is based on a thinning of the vertices in $g_{r}(\XS)$. To show the reverse direction, consider boxes $Q_s(sz)$ for $z\in \Z^d$ with $s>r$ and define for any set $A\subset\R^d$ the diameter of $A$ by $\diam(A)=\sup\{|x-y|\colon x,y\in A\}$. We say that a site $z\in \Z^d$ is a {\em good site} if all the following events happen
\begin{align*}
A_m(z)&=\{\text{for all }X_i\in \XS\cap Q_{3s}(sz)\colon \#\big(\XW\cap B_r(X_i)\big)\le m\},\\
B_n(z)&=\{\text{for all }X_i\in \XS\cap Q_{3s}(sz)\colon \#\big(\XS\cap B_r(X_i)\big)\le n+1\},\\
C_{\laI}(z)&=\{\text{all }X_i\in \XS\cap Q_{3s}(sz)\text{ are good nodes}\},\\
D_s(z)&=\{g_r(\XS)\text{ contains a unique cluster }Z\text{ in }Q_s(sz)\text{ with }\diam(Z)\ge s/2\}, \text{ and }\\
E_s(z)&=\{g_r(\XS)\text{ contains a unique cluster }Z\text{ in }Q_{3s}(sz)\text{ with }\diam(Z)\ge s/2\}.
\end{align*}
Otherwise $z$ is called a {\em bad site}. In particular, if $o$ is connected in $g_r(\XS)$ to $\partial Q_s$ and $o$ is contained in an infinite cluster of good sites $z\in \Z^d$, then the infection survives globally. 
Indeed, since $o$ is connected in $g_r(\XS)$ to $\partial Q_s$, there is a cluster $Z_o$ with $\diam(Z_o)\ge s/2$ in $Q_s(o)$. Further, let $z$ be a neighbor of $o$ in the infinite component of good sites. Then, there exists a unique cluster $Z_z$ in $Q_s(z)$ again with $\diam(Z_z)\ge s/2$. Since $Z_o$ and $Z_z$ are also unique in $Q_{3s}(z)$ respectively $Q_{3s}(o)$, $Z_o$ and $Z_z$ must be connected in $Q_{3s}(o)\cup Q_{3s}(z)$. This can be iterated along the path of good sites to infinity. By the goodness of that path, also there is no thinning, and hence the infection can globally survive.
Now we can estimate,
\begin{align*}
0<\theta(\muS)\le \P(o &\text{ is part of a finite cluster of good sites})+ \P(G),
\end{align*}
and it suffices to show that the percolation probability for the process of good sites can be pushed arbitrarily close to one as the parameters $s,m,n$ and $\laI$ tend to infinity. 
Note that by the definition of goodness of nodes, goodness of sites, and since $s>r$, the process of good sites is a $3$-dependent percolation process. 
Using the domination-by-product measure result~\cite[Theorem 0.0]{domProd}, it suffices to bound the $3$-dependent percolation process from below by a supercritical Bernoulli percolation process with parameter arbitrarily close to one, in the usual sense of evaluations of increasing events. In other words, it suffices to show that
\begin{align*}
\limsup_{s\uparrow\infty}&\limsup_{m\uparrow\infty}\limsup_{n\uparrow\infty}\limsup_{\laI\uparrow\infty}\P(o \text{ is a bad site})=0.
\end{align*}
For this, we can bound the probability for a bad site by
\begin{align*}
\P(o \text{ is a bad site})&= 1-\P(A_m\cap B_n\cap C_{\laI}\cap D_s\cap E_s)\\
&\le \P(A_m\cap B_n\cap C_{\laI}^{\rm c})+ \P(A_m^{\rm c})+\P(B_n^{\rm c})+\P(D_s^{\rm c})+\P(E_s^{\rm c}),
\end{align*}
where we suppressed the dependence on $z=o$. 
Now, using the large-deviation estimates in~\cite[Theorem 2]{uniqFin}, we can choose $s$ sufficiently large such that 
\begin{align*}
\P(D_s^{\rm c})<\e\qquad \text{ and }\qquad\P(E_s^{\rm c})<\e. 
\end{align*}
Next, for given $s$, we can choose $n$ and $m$ sufficiently large such that also
\begin{align*}
\P(A_m^{\rm c})<\e\qquad \text{ and }\qquad\P(B_n^{\rm c})<\e,
\end{align*}
by convergence in bounded domains. Finally, for given $s,m$ and $n$, note that, under the events $A_m$ and $B_n$, the probability of a node $X_i\in \XS\cap Q_{3s}(o)$ to be an open node is bounded from below by
\begin{align*}
\left(\frac{\laI}{\laI+\#\big(B_r(X_i) \cap \XS\big)-1+\#\big(B_r(X_i) \cap \XW\big)}\right)^{\#(B_r(X_i) \cap \XS)-1}\ge \left(\frac{\laI}{\laI+n+m}\right)^n. 
\end{align*}
Moreover, by the neighbor constraint imposed by the event $B_n$, there is also a maximal number of nodes that can be contained in $Q_{3s}(o)$, i.e., there exists $k=k(s,n)\in\N$ such that $\#(\XS\cap Q_{3s}(o))<k$. This implies that 
$$\P(A_m\cap B_n\cap C_{\laI}^{\rm c})\le1-\left(\frac{\laI}{\laI+n+m}\right)^{nk},$$
where we used that the indicators that nodes are open is a family of independent random variable indexed by the nodes in $\XS$. In particular, for given $s,m$ and $n$, we can now choose $\laI$ sufficiently large such that also 
$$\P(A_m\cap B_n\cap C_{\laI}^{\rm c})<\e,$$
which concludes the proof. 
\end{proof}


\bibliography{../AlexBene.bib}

\begin{thebibliography}{MMVY16}

\bibitem[And92]{An92}
E.D. Andjel.
\newblock Survival of multidimensional contact process in random environments.
\newblock {\em Bol. Soc. Brasil. Mat. (N.S.)}, 23(1-2):109--119, 1992.

\bibitem[CD09]{ChDu09}
S.~Chatterjee and R.~Durrett.
\newblock Contact processes on random graphs with power law degree
  distributions have critical value 0.
\newblock {\em Ann. Probab.}, 37(6):2332--2356, 2009.

\bibitem[DJT18]{DuJuTa18}
R.~Durrett, M.~Junge, and S.~Tang.
\newblock Coexistence in chase-escape.
\newblock {\em arXiv preprint arXiv:1807.05594}, 2018.

\bibitem[Dur10a]{Du10a}
R.~Durrett.
\newblock {\em Random graph dynamics}, volume~20 of {\em Cambridge Series in
  Statistical and Probabilistic Mathematics}.
\newblock Cambridge University Press, Cambridge, 2010.

\bibitem[Dur10b]{Du10b}
R.~Durrett.
\newblock Some features of the spread of epidemics and information on a random
  graph.
\newblock {\em Proceedings of the National Academy of Sciences},
  107(10):4491--4498, 2010.

\bibitem[Kor05]{Ko05}
G.~Kordzakhia.
\newblock The escape model on a homogeneous tree.
\newblock {\em Electron. Commun. Probab.}, 10:113--124, 2005.

\bibitem[Lig85]{Li85}
T.M. Liggett.
\newblock {\em Interacting Particle Systems}.
\newblock Grundlehren der mathematischen Wissenschaften. Springer-Verlag, 1985.

\bibitem[Lig92]{Li92}
T.M. Liggett.
\newblock The survival of one-dimensional contact processes in random
  environments.
\newblock {\em Ann. Probab.}, 20(2):696--723, 1992.

\bibitem[LSS97]{domProd}
T.M. Liggett, R.H. Schonmann, and A.M. Stacey.
\newblock Domination by product measures.
\newblock {\em Ann. Probab.}, 25(1):71--95, 1997.

\bibitem[MMVY16]{MoMoVaYa16}
T.~Mountford, J.-C. Mourrat, D.~Valesin, and Q.~Yao.
\newblock Exponential extinction time of the contact process on finite graphs.
\newblock {\em Stochastic Process. Appl.}, 126(7):1974--2013, 2016.

\bibitem[MS16]{MeSi16}
L.~M\'{e}nard and A.~Singh.
\newblock Percolation by cumulative merging and phase transition for the
  contact process on random graphs.
\newblock {\em Ann. Sci. \'{E}c. Norm. Sup\'{e}r. (4)}, 49(5):1189--1238, 2016.

\bibitem[MV16]{MoVa16}
J.-C. Mourrat and D.~Valesin.
\newblock Phase transition of the contact process on random regular graphs.
\newblock {\em Electron. J. Probab.}, 21:Paper No. 31, 17, 2016.

\bibitem[MVY13]{MoVaYa13}
T.~Mountford, D.~Valesin, and Q.~Yao.
\newblock Metastable densities for the contact process on power law random
  graphs.
\newblock {\em Electron. J. Probab.}, 18:No. 103, 36, 2013.

\bibitem[Pen91]{pcPerc}
M.D. Penrose.
\newblock On a continuum percolation model.
\newblock {\em Adv. in Appl. Probab.}, 23(3):536--556, 1991.

\bibitem[PP96]{uniqFin}
M.D. Penrose and A.~Pisztora.
\newblock Large deviations for discrete and continuous percolation.
\newblock {\em Adv. in Appl. Probab.}, 28(1):29--52, 1996.

\bibitem[TKL18]{TaKoLa18}
S.~Tang, G.~Kordzakhia, and S.P. Lalley.
\newblock Phase transition for the chase-escape model on 2d lattices.
\newblock {\em arXiv preprint arXiv:1807.08387}, 2018.

\end{thebibliography}
\bibliographystyle{alpha}

\end{document}